\documentclass[10pt]{article}
\usepackage{float}
\usepackage{graphicx}
\usepackage{epstopdf}
\usepackage{cite}
\usepackage{mathrsfs}
\usepackage{amsfonts}
\usepackage{amsmath}
\usepackage{amsfonts,amssymb,color}
\usepackage{dsfont}
\usepackage{curves}
\usepackage{mathrsfs}
\usepackage{pifont}
\usepackage{amssymb}
\usepackage{tikz}
\usepackage{hyperref}
\usepackage{enumitem}
\usepackage{CJK}
\usepackage[english]{babel}
\usepackage{times}
\usepackage[T1]{fontenc}
\usepackage{indentfirst}
\usepackage{subfigure}

\usepackage{subfigure}
\usepackage[T1]{fontenc}
\usepackage[utf8]{inputenc}
\usepackage{authblk}
\usepackage[misc]{ifsym}
\newtheorem{theorem}{\bf Theorem}[section]

\newtheorem{corollary}[theorem]{\bf Corollary}

\newcommand{\qed}{\hfill\rule{0.5em}{0.809em}}

\def\pf{\noindent {\it Proof. }}

\def\qed{\hfill \rule{4pt}{7pt}}

\date{}
\begin{document}

\title{ New sufficient degree conditions for an $r$-uniform hypergraph to be $k$-edge-connected
 \thanks{Supported by National Natural
        Science Foundation of China (Grant No. 11961019) and Hainan Provincial Natural Science Foundation
        of China (Grant No. 621RC510).} }
\author[a,b]{Jiyun Guo}
\author[a]{Jun Wang}
\author[a]{Zhanyuan Cai}
\author[b]{Haiyan Li \footnote{Corresponding author: lhy9694@163.com}}
\affil[a]{Center for Applied Mathematics, Tianjin University,
Tianjin, 300072, China} \affil[b]{School of Science, Hainan
University, Haikou, Hainan, 570228, China}
    \maketitle

 \begin{abstract}

An $r$-uniform hypergraphic sequence (i.e., $r$-graphic sequence)
$d=(d_1, d_2,\cdots,d_n)$ is said to be forcibly $k$-edge-connected
if every realization of $d$ is $k$-edge-connected. In this paper, we
obtain a strongest sufficient degree condition for $d$ to be
$k$-edge-connected for all $k\ge 1$ and a strongest sufficient
degree condition for $d$ to be super edge-connected. As a corollary,
we give the minimum degree condition for $d$ to be maximally
edge-connected. We also obtain another sufficient degree condition
for $d$ to be $k$-edge-connected.

 \begin{flushleft}
 {\em Keywords and phrases:} $r$-uniform hypergraph; $k$-edge-connected;
 super edge-connected; monotone degree condition

{\em Mathematics Subject Classification:} 05C07
 \end{flushleft}

\end{abstract}

\section{Introduction}
First let us introduce some terminology and notations. Unless other
stated, we follow [2] for undefined terms on hypergraphs. Let
$H=(V,E)$ be a hypergraph on $n$ vertices with vertex set $V$ and
edge set $E$. We say that $H$ is an $r$-uniform hypergraph if every
edge contains $r$ vertices, and that $H$ is a complete $r$-uniform
hypergraph, denoted by $K_n^r$, if $E$ consists of all $r$-subsets
of $V$. $H$ is simple if there are no repeated edges. Thus, a simple
2-uniform hypergraph is a simple graph. For a vertex $v$ in
hypergraph $H$, the degree of $v$, denoted $d_H(v)$ (or simply
$d(v)$ when $H$ is understood) is the number edges of $H$ that
contain $v$.

Let $d=(d_1, d_2, \cdots, d_n)$ be a sequence of nonnegative
integers with $d_1\le d_2\le \cdots \le d_n$. We say that $d$ is
$r$-uniform hypergraphic if it is the degree sequence of a simple
$r$-uniform hypergraph $H$ on $n$ vertices, and such an $r$-uniform
hypergraph $H$ is referred to as a realization of $d$. When $r=2$,
we will simply say that $d$ is graphic. The degree sequences for
simple graphs have been studied for many years and and by several
authors, including the celebrated work of Erd\"{o}s and Gallai [8].
Based on this, Sierksma and Hoogeveen [12] listed seven criteria and
 Cai et al.[5] gave it a generalization.

In this article, we mainly study sufficient degree conditions for a
simple hypergraph to be $k$-edge-connected and super edge-connected.
We say that a hypergraph $H=(V,E)$ is connected if for any pair of
vertices $u, v\in V$ there exists a path from $u$ to $v$; otherwise
we say $H$ is disconnected. For any $E_0\subseteq E$, if $H-E_0=H
(V, E-E_0)$ is disconnected, we say that $E_0$ is an edge-cut. The
edge connectivity of $H$, denoted by $\lambda (H)$, is the least
cardinality of all edge-cuts of $H$. A hypergraph $H$ is
$k$-edge-connected if $\lambda(H)\ge k$, where $1\le k\le \delta(H)$
with $\delta(H)$ denoting the minimum degree of $H$. If
$\lambda(H)=\delta(H)$, we say that $H$ is maximally edge-connected.
Moreover, a hypergraph $H$ is super edge-connected (or,
super-$\lambda$) if every minimum edge-cut consists of edges
incident with one vertex with minimum degree. There are several
degree sequence conditions for maximally edge-connected graphs and
super edge-connected graphs. For example, Dankelmann and Meierling
[7] and Zhao et al. [15], etc.

If $d=(d_1, d_2,\cdots,d_n)$ and $d'=(d'_1, d'_2,\cdots,d'_n)$ are
two sequences of nonnegative integers, we say that $d'$ majorizes
$d$, denoted $d'\ge d$, if $d'_j\ge d_j$ for $1\le j\le n$. Let $P$
denote an $r$-uniform hypergraphic property and $d$ be an
$r$-uniform hypergraphic sequence. If every realization  of $d$ has
property $P$, we say that $d$ is forcibly $P$. Historically,
sufficient conditions for a graphic sequence to have a certain
property, e.g., $k$-connected [3,4], $k$-edge-connected [1,10,14],
hamiltonian [6], clique size [13], etc, have been studied. We call a
forcibly $P$ degree condition $\chi$ is monotone increasing if
$d',d$ are hypergraphic, $d'\ge d$ and $d$ satisfies $\chi$ implies
that $d'$ satisfies $\chi$. Furthermore, $\chi$ is said to be a
strongest degree condition for property $P$ if whenever it does not
guarantee that $d$ is forcibly $P$, then $d$ is majorized by an
$r$-uniform hypergraphic sequence $d'$ which has a realization
without property $P$. In 1969, Bondy [4] gave a strongest degree
condition for a graphic sequence to be forcibly $k$-connected. In
1972, Chv$\acute{a}$tal [6] gave a strongest degree condition for a
graphic sequence to be forcibly hamiltonian.  Recently, Frosini et
al. [9] gave new sufficient conditions on the degree sequences of
uniform hypergraphs. Liu et al. [11] gave a simple sufficient degree
conditions for a uniform hypergraph to be $k$-edge-connected or
super edge-connected and the strongest monotone increasing degree
conditions for a uniform hypergraph to be $k$-edge-connected when
$k=1, 2, 3$. Inspired by this, we give the strongest sufficient
degree conditions for an $r$-uniform hypergraph to be
$k$-edge-connected for all $k\ge 1$ and super edge-connected. These
are in section 2 and 3. Additionally, Yin et al. gave a sufficient
degree condition for a graphic sequence to be forcibly
$k$-edge-connected. We shall extend the result to hypergraphs in
section 4.


\section{The strongest conditions for $k$-edge-connected hypergraphs}

The following theorem due to Liu, Meng and Tian [11] gave simple
sufficient degree conditions for an $r$-uniform hypergraph to be
$k$-edge-connected. To present the result, they utilized a parameter
$g$ and define $g=g(k, r)(\ge r)$ to be the smallest integer such
that $gk\le g {g-1\choose r-1}+(r-1)(k-1)$, i.e., $k\le {g\choose
r-1}-{r-1 \choose g-r+1}$, where $r\ge 2$ and $k\ge 2$.

 \begin{theorem} ([11])\label{21t}
For two integers $r\ge2$ and $k\ge 2$. Let $d=(d_1, d_2, \cdots,
d_n)$ be an $r$-uniform hypergraphic sequence with $d_1\le
d_2\le\cdots\le d_n$. If $d$ satisfies all of the following
conditions, then $d$ is forcibly $k$-edge connected:

\begin{itemize}
\item [$(1)$]  $d_1\ge k$,
\item [$(2)$] $d_{j-(r-1)(k-1)}\le {j-1\choose r-1}$ and $d_j\le {j-1\choose r-1}+k-1$ implies $d_n\ge{n-j-1\choose r-1}+k$ for $g(k, r)\le j\le \lfloor {n\over 2}\rfloor$.
\end{itemize}
\end{theorem}

In addition to Theorem \ref{21t}, Liu et al. also established the
strongest monotone increasing conditions for an $r$-uniform
hypergraph be be $k$-edge-connected when $k=1, 2, 3$. Now comes the
case for $k=2$, which differs a little from one of the corollaries
of Theorem 2.3. In fact, however, they are equivalent.

 \begin{theorem} ([11])\label{22t}
Given integer $r\ge2$. Let $d=(d_1, d_2, \cdots, d_n)$ be an
$r$-uniform sequence with $d_1\le d_2\le\cdots\le d_n$. If $d$
satisfies all of the following conditions, then $d$ is forcibly
$2$-edge connected:

\begin{itemize}
\item [$(1)$]  $d_1\ge 2$,
\item [$(2)$] $d_{j-t_1}\le {j-1\choose r-1}$ and $d_j\le {j-1\choose r-1}+1$ implies
 $d_{n-r+t_1}\ge{n-j-1\choose r-1}+1$ or $d_{n}\ge{n-j-1\choose r-1}+2$
 for $r+1\le j< {n\over 2}$ and $1\le t_1\le r-1$.

\item [$(3)$] $d_{n\over 2}\le {{n\over 2}-1\choose r-1}$ and $d_{n-r}\le {{n\over 2}-1\choose r-1}$
implies $d_{n}\ge{{n\over 2}-1\choose r-1}+2$ for even $n\ge 2r+2$.

\end{itemize}
\end{theorem}

Inspired by these results, we first give  the strongest monotone
degree condition
 for an $r$-uniform hypergraph to be $k$-edge-connected for all $k\ge 1$, in which four conditions are contained, where (2)-(4) are "Bondy-Chv$\acute{a}$tal type" conditions.
Before
 presenting it, we introduce some parameters.

Given integers $r\ge 2$ and $k\ge 2$, $j^\ast=j^\ast (k,r)(r+1\le
j^\ast <{n\over 2})$ is defined to be the biggest integer such that
${j^\ast-1\choose r-1}+k-1\le {n-j^\ast-1\choose r-1}$. Let
$t_1,\cdots,t_{k-1},s_1,\cdots,$  $s_{k-1}$ be nonnegative integers
which satisfied:
$$[t_1+2t_2+\cdots+(k-1)t_{k-1}]+[s_1+2s_2+\cdots+(k-1)s_{k-1}]=(k-1)r,$$
 $$k-1\le t_1+s_1\le (k-1)r,$$
 $$0\le t_i+s_i\le \lfloor{k-1\over i}\rfloor(r-1),2\le i \le k-1,$$
$$1\le t_1+t_2+\cdots+t_{k-1}\le \min\{(k-1)r-(k-1), j\}, k\ge
2,$$
$$0\le t_i+\cdots+t_{k-1}\le \lfloor{k-1\over i}\rfloor(r-1), k\ge 3,2\le i\le k-2,$$
$$1\le s_1+s_2+\cdots+s_{k-1}\le \min\{(k-1)r-(k-1), n-j\}, k\ge 2,$$
$$0\le s_i+\cdots+s_{k-1}\le \lfloor{k-1\over i}\rfloor(r-1), k\ge 3,2\le i\le k-2.$$

 \begin{theorem} \label{23t}
    For integers $r\ge2$ and $k\ge 1$. Let $d=(d_1, d_2, \cdots, d_n)$ be an
    $r$-uniform hypergraphic sequence with $d_1\le d_2\le\cdots\le d_n$.
    If $d$ satisfies all of the following conditions, then $d$ is forcibly $k$-edge-connected:

    \begin{itemize}
        \item [$(1)$]  $d_1\ge k$,
        \item [$(2)$] $d_{j-t_1-t_2-\cdots-t_{k-1}}\le {j-1\choose r-1}$, $d_{j-t_2-\cdots-t_{k-1}}
        \le {j-1\choose r-1}+1$,
        $\cdots,$ $d_{j-t_{k-1}}\le {j-1\choose r-1}+k-2$ and $d_j\le {j-1\choose r-1}+k-1$
        implies $d_{n-s_1-\cdots-s_{k-1}}\ge{n-j-1\choose r-1}+1$ or
        $d_{n-s_2-\cdots-s_{k-1}}\ge{n-j-1\choose r-1}+2$, $\cdots$, or
        $d_{n-s_{k-1}}\ge{n-j-1\choose r-1}+k-1$ or $d_{n}\ge{n-j-1\choose r-1}+k$
         for $r+1\le j\le j^\ast$,
\item [$(3)$] $d_{j-t_1-\cdots-t_{k-1}}\le
         {n-j-1\choose r-1}-\vartriangle$, $d_{j-t_2-\cdots-t_{k-1}}
         \le {n-j-1\choose r-1}+1-\vartriangle$,$\cdots$,
         $d_{j-t_\vartriangle-\cdots-t_{k-1}} \le {n-j-1\choose r-1}-1$,
         $\cdots$, $d_{j-t_{k-1}} \le {n-j-1\choose
         r-1}+k-2-\vartriangle$ and $d_j \le {n-j-1\choose
         r-1}+k-1-\vartriangle$ implies
          $d_{n-t_{\vartriangle+1}-\cdots-t_{k-1}-s_1-\cdots-s_{k-1}} \ge {n-j-1\choose
          r-1}+1$, $\cdots$, or
          $d_{n-t_{k-1}-s_{k-1-\vartriangle}-\cdots-s_{k-1}} \ge {n-j-1\choose r-1}+k-1-\vartriangle$
or
          $d_{n-s_{k-\vartriangle}-\cdots-s_{k-1}} \ge {n-j-1\choose
          r-1}+k-\vartriangle $, $\cdots$, or $d_{n-s_{k-1}} \ge {n-j-1\choose
          r-1}+k-1$ or $d_n \ge {n-j-1\choose
          r-1}+k$  for  $j^\ast<j<{n\over 2}$,
          where $\vartriangle={n-j-1\choose  r-1}-{j-1\choose
          r-1}$.

        \item [$(4)$] $d_{n-t_1-\cdots-t_{k-1}-s_1-\cdots-s_{k-1}}\le
         {{n\over 2}-1\choose r-1}$, $d_{n-t_2-\cdots-t_{k-1}-s_2-\cdots-s_{k-1}}
         \le {{n\over 2}-1\choose r-1}+1$,$\cdots$, and $d_{n-t_{k-1}-s_{k-1}}
         \le {{n\over 2}-1\choose r-1}+k-2$ implies $d_n
         \ge {{n\over 2}-1\choose r-1}+k$ for even $n\ge 2r+2$.

    \end{itemize}
\end{theorem}

\pf By induction on $k$. The theorem is trivial for $k=1, 2, 3$
(refer to Corollary 2.4-2.6). Assume the theorem holds for all $k\le
m-1$. Now suppose that $d$ satisfies $(1)-(4)$ for $k=m\ge 4$, but
$d$ is not forcibly $m$-edge-connected. Namely, there exists a
non-$m$-edge-connected $r$-uniform realization $H$ consisting of two
subhypergraphs $C_1$ and $C_2$ joined by $m-1$ hyperedges $e_1, e_2,
\cdots, e_{m-1}$, where $V(C_2)=V(H)\setminus V(C_1)$ and
$|V(C_1)|\le |V(C_2)|$. Denote $E_0=\{e_1, e_2, \cdots, e_{m-1}\}$
and let $|V(C_1)|=j$. One can see that $r+1\le j\le \lfloor {n\over
2}\rfloor$, since $j\le r$ implies there exists a vertex of degree
less than $m$. Let $U_i(V_i)$ be the set of vertices in
$V(H_1)(V(H_2))$, any one of which belong to exactly $t_i(s_i)$
hyperedges of $e_1, e_2,\cdots, e_{m-1}$, where $1\le i\le m-1$.
Then we can see that
$[t_1+2t_2+\cdots+(m-1)t_{m-1}]+[s_1+2s_2+\cdots+(m-1)s_{m-1}]=(m-1)r$,
$m-1\le t_1+s_1\le (m-1)r$, $0\le t_i+s_i\le \lfloor{m\over
i}\rfloor(r-1)$ for $i=2,\cdots,m-1$, $1\le
t_1+t_2+\cdots+t_{m-1}\le \min\{(m-1)r-(m-1), j\}$,  $0\le
t_2+\cdots+t_{m-1}\le \lfloor{m-1\over 2}\rfloor(r-1),\cdots,0\le
t_{m-2}+t_{m-1}\le \lfloor{m-1\over m-2}\rfloor(r-1)$, $1\le
s_1+s_2+\cdots+s_{m-1}\le \min\{(m-1)r-(m-1), n-j\}$, $0\le
s_2+\cdots+s_{m-1}\le \lfloor{m-1\over 2}\rfloor(r-1)$, $\cdots$,
$0\le s_{m-2}+s_{m-1}\le \lfloor{m-1\over m-2}\rfloor(r-1)$.

It is not difficult to see that there are at most
$t_1+t_2+\cdots+t_{m-1}$ vertices in $C_1$ of degree larger than
$j-1 \choose r-1$, at most $t_2+\cdots+t_{m-1}$ vertices in $C_1$ of
degree larger than ${j-1 \choose r-1}+1$, at most
$t_3+\cdots+t_{m-1}$ vertices in $C_1$ of degree larger than ${j-1
\choose r-1}+2$, $\cdots$, and at most  $t_{m-1}$ vertices in $C_1$
of degree larger than ${j-1 \choose r-1}+m-2$, while the degree of
each vertex in $C_1$ is at most ${j-1 \choose r-1}+m-1$. Thus we
have $d_{j-t_1-t_2-\cdots-t_{m-1}}\le {j-1\choose r-1}$ if $j>
t_1+t_2+\cdots+t_{m-1}$, $d_{j-t_2-\cdots-t_{m-1}}\le {j-1\choose
r-1}+1$ if $j> t_2+\cdots+t_{m-1}$, $\cdots$, $d_{j-t_{m-1}}\le
{j-1\choose r-1}+m-2$ if $j>t_{m-1}$ and $d_{j}\le {j-1\choose
r-1}+m-1$. We distinguish the following cases.

{\bf Case  1.}  $r+1\le j<{n\over 2}$.

{\bf Case  1.1.} $r+1\le j\le j^\ast .$

 In this case, each vertex in
$C_1$ has degree at most ${n-j-1\choose r-1}$
 according to  the definition of $j^\ast$ described above.
Additionally, there are at most $s_1+ \cdots + s_{m-1}$ vertices in
$C_2$ of degree larger than ${n-j-1\choose r-1}$, at most  $s_2+
\cdots + s_{m-1}$ vertices in $C_2$ of degree larger than
${n-j-1\choose r-1}+1, \cdots,$ and at most  $s_{m-1}$ vertices in
$C_2$ of degree larger than ${n-j-1\choose r-1}+m-2$. Thus,
$d_{n-s_1-\cdots-s_{m-1}}\le{n-j-1\choose r-1}$,
$d_{n-s_2-\cdots-s_{m-1}}\le{n-j-1\choose r-1}+1$, $\cdots$,
$d_{n-s_{m-1}}\le{n-j-1\choose r-1}+m-2$ and $d_{n}\le{n-j-1\choose
r-1}+m-1$. Therefore the consequent in $(2)$ fails, while the
antecedent in $(2)$ is satisfied, yielding a contradiction.

{\bf Case  1.2.} $j^\ast<j <{n\over 2}.$

Obviously, ${j-1\choose r-1}<{n-j-1\choose r-1}$. Write
$\vartriangle={n-j-1\choose r-1}-{j-1\choose r-1}$, then $1\le
\vartriangle <m-1$ since ${j-1\choose r-1}+m-1>{n-j-1\choose r-1}$
when $j>j^\ast$. Note that the degree of each vertex in $V(C_1)-E_0$
is at most ${j-1\choose r-1}={n-j-1\choose r-1}-\vartriangle$, the
degree of each vertex in $ U_1$ is at most ${j-1\choose
r-1}+1={n-j-1\choose r-1}+1-\vartriangle$, the degree of each vertex
in $ U_2$ is at most ${j-1\choose r-1}+2={n-j-1\choose
r-1}+2-\vartriangle$, $\cdots$, the degree of each vertex in $
U_\vartriangle$ is at most ${j-1\choose
r-1}+\vartriangle={n-j-1\choose r-1}$, the degree of each vertex in
 $ U_{\vartriangle +1}$ is at most ${j-1\choose
r-1}+\vartriangle+1={n-j-1\choose r-1}+1$, $\cdots$, the degree of
each vertex in $U_{m-1}$ is at most ${j-1\choose
r-1}+m-1={n-j-1\choose r-1}+m-1-\vartriangle$, and so
$d_{j-t_1-\cdots-t_{m-1}}\le
         {n-j-1\choose r-1}-\vartriangle$, $d_{j-t_2-\cdots-t_{m-1}}
         \le {n-j-1\choose r-1}+1-\vartriangle$, $d_{j-t_3-\cdots-t_{m-1}}
         \le {n-j-1\choose r-1}+2-\vartriangle$, $\cdots$,
         $d_{j-t_\vartriangle-\cdots-t_{m-1}} \le {n-j-1\choose r-1}-1$,
         $d_{j-t_{\vartriangle+1}-\cdots-t_{m-1}} \le {n-j-1\choose
         r-1}$, $\cdots$, $d_{j-t_{m-1}} \le {n-j-1\choose
         r-1}+m-2-\vartriangle$ and $d_j \le {n-j-1\choose
         r-1}+m-1-\vartriangle$. Furthermore, the degree of each
         vertex in $V_i$ is at most ${n-j-1\choose
         r-1}+i$, where $1\le i\le m-1$. Thus,
         $d_{n-t_{\vartriangle+1}-\cdots-t_{m-1}-s_1-\cdots-s_{m-1}} \le {n-j-1\choose
          r-1}$, $d_{n-t_{\vartriangle+2}-\cdots-t_{m-1}-s_2-\cdots-s_{m-1}} \le {n-j-1\choose
          r-1}+1$, $\cdots$,
          $d_{n-t_{m-1}-s_{m-1-\vartriangle}-\cdots-s_{m-1}} \le {n-j-1\choose
          r-1}+m-2-\vartriangle$,
          $d_{n-s_{m-\vartriangle}-\cdots-s_{m-1}} \le {n-j-1\choose
          r-1}+m-1-\vartriangle $, $\cdots$,  $d_{n-s_{m-1}} \le {n-j-1\choose
          r-1}+m-2$ and $d_n \le {n-j-1\choose
          r-1}+m-1$, this contradicts (3).

{\bf Case  2.} $j={n\over 2}$ for even $n.$

 Note that the degree of each
vertex in $V(H)-E_0$ is at most ${n\over 2}-1 \choose r-1$.
Moreover, there are at most $t_2+\cdots+t_{m-1}+s_2+\cdots+s_{m-1}$
vertices in $H$ of degree larger than ${{n\over 2}-1 \choose r-1}+1,
\cdots$, and at most $s_{m-1}+t_{m-1}$ vertices in $H$ of degree
larger than ${{n\over 2}-1 \choose r-1}+m-2$, while the degree of
each vertex in $H$ is at most ${{n\over 2}-1 \choose r-1}+m-1$.
Thus, $d_{n-t_1-t_2-\cdots-t_{m-1}-s_1-\cdots-s_{m-1}}\le {{n\over
2}-1\choose r-1}$, $d_{n-t_2-\cdots-t_{m-1}-s_2-\cdots-s_{m-1}}\le
{{n\over 2}-1\choose r-1}+1, \cdots$, $d_{n-t_{m-1}-s_{m-1}}\le
{{n\over 2}-1\choose r-1}+m-2$ and $d_n \le {{n\over 2}-1\choose
r-1}+m-1$, contrary to $(4)$.

Therefore, the hypothesis is not valid, and so Theorem 2.3 follows
by the principle of induction. \qed

If $k=1, 2, 3$ in Theorem \ref{23t},
 then we can obtain the strongest monotone increasing degree conditions for
  an $r$-uniform hypergraph to be $k$-edge-connected when $k=1, 2, 3$.

 \begin{corollary}\label{24c}([11])
For integers $r\ge2$ and $k\ge 1$. Let $d=(d_1, d_2, \cdots, d_n)$
be an $r$-uniform hypergraphic sequence with $d_1\le d_2\le\cdots\le
d_n$. If $d_1\ge 1$ and $d_j \le {j-1 \choose r-1}$ implies $d_n \ge
{n-j-1 \choose r-1}+1$ for $r+1\le j\le \lfloor {n\over 2}\rfloor$,
then $d$ is forcibly connected.
\end{corollary}

 \begin{corollary}\label{25c}
    For integers $r\ge2$. Let $d=(d_1, d_2, \cdots, d_n)$ be an $r$-uniform
    hypergraphic sequence with $d_1\le d_2\le\cdots\le d_n$. If $d$ satisfies all of the
    following conditions, then $d$ is forcibly $2$-edge-connected:

    \begin{itemize}
        \item [$(1)$]  $d_1\ge 2$,
        \item [$(2)$] $d_{j-t_1}\le {j-1\choose r-1}$ and $d_j\le {j-1\choose r-1}+1$ implies $d_{n-s_1}\ge{n-j-1\choose r-1}+1$ or $d_{n}\ge{n-j-1\choose r-1}+2$ for $r+1\le j< {n\over 2}$, $1\le t_1\le r-1$ and $t_1+s_1=r$.

        \item [$(3)$] $d_{n-r}\le {{n\over 2}-1\choose r-1}$ implies $d_{n}\ge{{n\over 2}-1\choose r-1}+2$ for even $n\ge 2r+2$.
\end{itemize}

\end{corollary}

Note that ${n\over 2}< n-r$ since $n\ge 2r+2$, then $d_{n\over 2}\le
d_{n-r}$. Therefore, Corollary \ref{25c} and Theorem \ref{22t} are
equivalent.

 \begin{corollary}\label{26c}([11])
    For integers $r\ge3$. Let $d=(d_1, d_2, \cdots, d_n)$ be an $r$-uniform hypergraphic sequence with $d_1\le d_2\le\cdots\le d_n$. If $d$ satisfies all of the following conditions, then $d$ is forcibly $3$-edge-connected:

    \begin{itemize}
        \item [$(1)$]  $d_1\ge 3$,
        \item [$(2)$] $d_{j-t_1-t_2}\le {j-1\choose r-1}$, $d_{j-t_2}\le {j-1\choose r-1}+1$ and
         $d_j\le {j-1\choose r-1}+2$ implies $d_{n-s_1-s_2}\ge{n-j-1\choose r-1}+1$ or
         $d_{n-s_2}\ge{n-j-1\choose r-1}+2$ or $d_{n}\ge{n-j-1\choose r-1}+3$ for
         $r+1\le j< {n\over 2}$, $2\le t_1+s_1\le 2r$, $0\le t_2+s_2\le r-1$,
         $1\le t_1+t_2\le min\{2r-2, j\}$, $1\le s_1+s_2\le min\{2r-2, n-j\}$,
         $t_1+2t_2+s_1+2s_2=2r$ and $t_i, s_i\ge 0$ for $1\le i\le
         2$,

        \item [$(3)$] $d_{n-t_1-t_2-s_1-s_2}\le {{n\over 2}-1\choose r-1}$, $d_{n-t_2-s_2}\le
         {{n\over 2}-1\choose r-1}+1$ implies $d_{n}\ge{{n\over 2}-1\choose r-1}+3$ for even $n\ge 2r+2$, where $t_i$ and $s_i (1\le i\le 2)$ are defined as shown in $(2)$.
    \end{itemize}

\end{corollary}

 \begin{theorem}\label{27b}
The degree conditions in Theorem \ref{23t} are strongest.
\end{theorem}

\pf In order to very that the degree conditions in Theorem 2.3 are
strongest, we now show that if any of (1)-(4) fails for $d$, then
$d$ is majorized by an $r$-uniform hypergraphic sequence $d'$ having
a realization $H'$ which is not $k$-edge-connected.

 Let $H_{t_1, \cdots, t_{k-1}, s_1, \cdots, s_{k-1}}(j, n-j)$
denote an $r$-uniform hypergraph consisting of two disjoint $K_j^r$
and $K_{n-j}^r$ joined by $k-1$ edges $e_1, e_2, \cdots, e_{k-1}$,
where $t_i$ and $s_i$ $(1\le i\le k-1)$ are defined as shown in
Theorem 2.3.

If $(1)$ fails, we take $d^\prime=d$ and $H^\prime=H$, where $H$ is
a non-$k$-edge-connected $r$-uniform realization of $d$ with minimum
degree less than $k$.

 If $(2)$ fails for $t_1, \cdots, t_{k-1}, s_1, \cdots, s_{k-1}$ and $j$ with $r+1\le j <j^\ast$,
 we take $d^\prime=({j-1 \choose r-1}^{j-t_1-\cdots-t_{k-1}}\big({j-1\choose r-1}+
 1\big)^{t_1}\cdots \big({j-1\choose r-1}
 +(k-1)\big)^{t_{k-1}}
 {n-j-1 \choose r-1}^{n-j-s_1-\cdots-s_{k-1}}
 \big({n-j-1\choose r-1}
 +1\big)^{s_1}\big({n-j-1\choose r-1}+2\big)^{s_2}\cdots \big({n-j-1\choose r-1}
 +(k-2)\big)^{s_{k-2}}\big({n-j-1\choose r-1}+(k-1)\big)^{s_{k-1}})$ and
 $H^\prime=H_{t_1, \cdots, t_{k-1}, s_1,\cdots,s_{k-1}}(j, n-j)$,
 where the symbol $x^y$ in $d'$ stands for $y$ consecutive terms,
 each equals to $x$.

If $(3)$ fails for $t_1,\cdots,t_{k-1}, s_1,\cdots,s_{k-1}$ and $j$
with $j^\ast<j<{n\over 2}$, we may take $d'=(({n-j-1 \choose
r-1}-\vartriangle)^{j-t_1-\cdots-t_{k-1}}({n-j-1 \choose
r-1}+1-\vartriangle)^{t_1}\cdots ({n-j-1 \choose
r-1}+k-1-\vartriangle)^{t_{k-1}} {n-j-1 \choose
r-1}^{n-j-t_{\vartriangle +1}-\cdots-t_{k-1}-s_1-\dots-s_{k-1}}
({n-j-1 \choose r-1}+1)^{t_{\vartriangle+1}}\cdots ({n-j-1 \choose
r-1}+k-1-\vartriangle)^{t_{k-1}}({n-j-1 \choose
r-1}+k-\vartriangle)^{s_{k-\vartriangle}}\cdots ({n-j-1 \choose
r-1}+k-1)^{s_{k-1}})$ and  $H^\prime=H_{t_1,\cdots,t_{k-1},
s_1,\cdots,s_{k-1}}$ $(j, n-j)$,

If $(4)$ fails for $t_1, \cdots, t_{k-1}, s_1, \cdots, s_{k-1}$, we
may take $d^\prime=({{n\over 2}-1 \choose
r-1}^{n-t_1-\cdots-t_{k-1}-s_1-\cdots-s_{k-1}} \big({{n\over
2}-1\choose r-1}+1\big)^{t_1+s_1} \big({{n\over 2}-1\choose
r-1}+2\big)^{t_2+s_2} \cdots \big({{n\over 2}-1\choose
r-1}+(k-2)\big)^{t_{k-2}+s_{k-2}} \big({{n\over 2}-1\choose
r-1}+(k-1)\big)^{t_{k-1}+s_{k-1}})$ and $H^\prime=H_{t_1, \cdots,
t_{k-1}, s_1, \cdots, s_{k-1}}$ $({n\over 2}, {n\over 2})$. \qed

\begin{corollary}\label{28c}
    For integers $r\ge 2$. Let $d=(d_1, d_2, \cdots, d_n)$ be an
    $r$-uniform hypergraphic sequence with $\delta=d_1\le d_2\le\cdots\le d_n$.
    If $d$ satisfies $(2)-(4)$ in Theorem \ref{23t} for $k=\delta$, then $d$ is forcibly maximally edge-connected.
\end{corollary}

Taking $t_1=(r-1)(k-1)$, $t_2=t_3=\cdots=t_{k-1}=0$, $s_1=
s_2=\cdots=s_{k-2}=0$ and $s_{k-1}=1$ in Theorem \ref{23t}, we
immediately get the following conclusion.

 \begin{corollary}\label{29c}
    For integers $r\ge 2$. Let $d=(d_1, d_2, \cdots, d_n)$ be an $r$-uniform hypergraphic sequence with $d_1\le d_2\le\cdots\le d_n$. If $d$ satisfies each of the following, then $d$ is forcibly $k$-edge-connected:

    \begin{itemize}
        \item [$(1)$]  $d_1\ge k$,
        \item [$(2)$] $d_{j-(r-1)(k-1)}\le {j-1\choose r-1}$ and $d_{j}\le {j-1\choose r-1}+k-1$ implies  $d_{n-1}\ge{n-j-1\choose r-1}+k-1$ or $d_{n}\ge{n-j-1\choose r-1}+k$ for $(r-1)(k-1)\le j \le \lfloor{n\over 2}\rfloor$.
    \end{itemize}
\end{corollary}

Recall that the parameter $g(k,r)\ge (r-1)(k-1)$ in Theorem
\ref{21t} since the index $j-(r-1)(k-1)>0$. Therefore, Corollary
\ref{29c} is a generalization of Theorem \ref{21t}.

\section{Super edge-connected $r$-uniform hypergraphs}

In this section, We provide the strongest degree conditions for an
$r$-uniform hypergraphic sequence with minimum degree $3$ or $k$ to
be forcibly super edge-connected.

\begin{theorem}\label{31t}
    For integer $r\ge 4$. Let $d=(d_1, d_2, \cdots, d_n)$ be an
    $r$-uniform hypergraphic sequence with $3=d_1\le d_2\le\cdots\le d_n$.
    If $d$ satisfies each of the following, then $d$ is forcibly super edge-connected:

    \begin{itemize}
    \item [$(1)$] $d_{j-t_1-t_2-t_3}\le {j-1\choose r-1}$,
    $d_{j-t_2-t_3}\le {j-1\choose r-1}+1$, $d_{j-t_3}\le {j-1\choose r-1}+2$ and
    $d_j\le {j-1\choose r-1}+3$ implies $d_{n-s_1-s_2-s_3}\ge{n-j-1\choose r-1}+1$
    or $d_{n-s_2-s_3}\ge{n-j-1\choose r-1}+2$ or $d_{n-s_3}\ge{n-j-1\choose r-1}+3$
    or $d_{n}\ge{n-j-1\choose r-1}+4$ for $r\le j< {n\over 2}$,

    \item [$(2)$] $d_{n-t_1-t_2-t_3-s_1-s_2-s_3}\le {{n\over 2}-1\choose r-1}$,
     $d_{n-t_2-t_3-s_2-s_3}\le {{n\over 2}-1\choose r-1}+1$,
     $d_{n-t_3-s_3}\le {{n\over 2}-1\choose r-1}+2$ implies $d_{n}\ge{{n\over 2}-1\choose r-1}+4$
     for even $n\ge 2r+2$,
 where $t_1+2t_2+3t_3+s_1+2s_2+3s_3=3r$,
    $3\le t_1+s_1\le 3r$,
     $0\le t_2+s_2\le 2(r-1)$, $0\le t_3+s_3\le r-1$,
     $1\le t_1+t_2+t_3\le min\{3r-3, j\}$,
     $1\le s_1+s_2+s_3\le min\{3r-3, n-j\}$,
     $0\le t_2+t_3\le r-1$, $0\le s_2+s_3\le r-1$ and
 $t_i(s_i)\ge 0$ for $1\le i\le 3$.
    \end{itemize}
\end{theorem}

\pf If $d$ satisfies $(1)$ and $(2)$, then any $r$-uniform
realization of $d$ is maximally edge-connected by Corollary
\ref{28c} corresponding to $k=3$. Now suppose that $d$ satisfies
$(1)$ and $(2)$ and there exists a non-super-$\lambda$ $r$-uniform
realization $H$ consisting of two subhypergraphs $C_1$ and $C_2$
joined by three edges $e_1, e_2$ and $e_3$, where
$V(C_2)=V(H)\setminus V(C_1)$ and $|V(C_1)|\le |V(C_2)|$. Denote
$E_0=\{e_1, e_2,e_3\}$ and let $|V(C_1)|=j$. Note that $r\le j \le
\lfloor{n\over 2}\rfloor$. Let $U_i(V_i)$ be the set of vertices in
$V(H_1)(V(H_2))$, any one of which belong to exactly $t_i(s_i)$
hyperedges of $e_1, e_2,\cdots, e_{m-1}$, where $1\le i\le m-1$.
Then we know that $t_1+2t_2+3t_3+s_1+2s_2+3s_3=3r$, $3\le t_1+s_1\le
3r$, $0\le t_2+s_2\le 2(r-1)$, $0\le t_3+s_3\le r-1$, $1\le
t_1+t_2+t_3\le min\{3r-3, j\}$, $1\le s_1+s_2+s_3\le min\{3r-3,
n-j\}$, $0\le t_2+t_3\le r-1$, $0\le s_2+s_3\le r-1$ and $t_i(
s_i)\ge 0$ for $1\le i\le 3$.

We can see that the degree of each vertex in $V(C_1)-E_0$ is at most
$j-1 \choose r-1$, the degree of each vertex in $U_1$ is at most
${j-1 \choose r-1}+1$, the degree of each vertex in $U_2$ is at most
${j-1 \choose r-1}+2$ and the degree of each vertex in $U_3$ is at
most ${j-1 \choose r-1}+3$. So we have $d_{j-t_1-t_2-t_3}\le
{j-1\choose r-1}$, $d_{j-t_2-t_3}\le {j-1\choose r-1}+1$,
$d_{j-t_3}\le {j-1\choose r-1}+2$ and $d_j\le {j-1\choose r-1}+3$.

If $j<{n\over 2}$, the degree of each vertex in $C_1$ is at most
${n-j-1\choose r-1}$ since ${j-1\choose r-1}+3\le {n-j-2\choose
r-1}+3\le {n-j-2\choose r-1}+{n-j-3\choose r-2}+2 \le {n-j-2\choose
r-1}+{n-j-3\choose r-2}+{n-j-3\choose r-3}={n-j-2\choose
r-1}+{n-j-2\choose r-2}={n-j-1\choose r-1}$ for $r\ge 4$. Moreover,
the degree of each vertex in $V(C_2)-E_0$ is at most ${n-j-1\choose
r-1}$, the degree of each vertex in $V_1$ is at most ${n-j-1\choose
r-1}+1$, the degree of each vertex in $V_2$ is at most
${n-j-1\choose r-1}+2$ and the degree of each vertex in $V_3$ is at
most ${n-j-1\choose r-1}+3$. Then we have
$d_{n-s_1-s_2-s_3}\le{n-j-1\choose r-1}$,
$d_{n-s_2-s_3}\le{n-j-1\choose r-1}+1$, $d_{n-s_3}\le{n-j-1\choose
r-1}+2$ and $d_{n}\le{n-j-1\choose r-1}+3$, contrary to $(1)$.

If $j={n\over 2}$ for even $n$, the degree of each vertex in
$V(H)-E_0$ is at most ${n\over 2}-1 \choose r-1$, the degree of each
vertex in $V_1\cup U_1$ is at most ${{n\over 2}-1 \choose r-1}+1$,
the degree of each vertex in $V_2\cup U_2$ is at most ${{n\over 2}-1
\choose r-1}+2$, and the degree of each vertex in $V_3\cup U_3$ is
at most ${{n\over 2}-1 \choose r-1}+3$. Therefore,
$d_{n-t_1-t_2-t_3-s_1-s_2-s_3}\le {{n\over 2}-1\choose r-1}$,
$d_{n-t_2-t_3-s_2-s_3}\le {{n\over 2}-1\choose r-1}+1$,
$d_{n-t_3-s_3}\le {{n\over 2}-1\choose r-1}+2$ and $d_{n}\le{{n\over
2}-1\choose r-1}+3$, contrary to $(2)$.\qed

\begin{theorem}\label{32t}
    For integers $r\ge 2$ and $\delta \ge 1$. Let $d=(d_1, d_2, \cdots, d_n)$ be an $r$-uniform hypergraphic sequence with $\delta=d_1\le d_2\le\cdots\le d_n$. If $d$ satisfies each of the following, then $d$ is forcibly super edge-connected:

    \begin{itemize}
\item [$(1)$] $d_{j-t_1-t_2-\cdots-t_{\delta}}\le {j-1\choose r-1}$, $d_{j-t_2-\cdots-t_{\delta}}
        \le {j-1\choose r-1}+1$,
        $\cdots$, $d_{j-t_{\delta}}\le {j-1\choose r-1}+\delta-1$ and $d_j\le {j-1\choose r-1}+\delta$
        implies $d_{n-s_1-\cdots-s_{\delta}}\ge{n-j-1\choose r-1}+1$ or
        $d_{n-s_2-\cdots-s_{\delta}}\ge{n-j-1\choose r-1}+2$,$\cdots$, or
        $d_{n-s_{\delta}}\ge{n-j-1\choose r-1}+\delta$ or $d_{n}\ge{n-j-1\choose r-1}+\delta+1$
         for $r+1\le j\le j^\ast$.
\item [$(2)$] $d_{j-t_1-\cdots-t_{\delta}}\le
         {n-j-1\choose r-1}-\vartriangle$, $d_{j-t_2-\cdots-t_{\delta}}
         \le {n-j-1\choose r-1}+1-\vartriangle$,$\cdots$,
         $d_{j-t_\vartriangle-\cdots-t_{\delta}} \le {n-j-1\choose r-1}-1$,
         $\cdots$, $d_{j-t_{\delta}} \le {n-j-1\choose
         r-1}+\delta-1-\vartriangle$ and $d_j \le {n-j-1\choose
         r-1}+\delta-\vartriangle$ implies
          $d_{n-t_{\vartriangle+1}-\cdots-t_{\delta}-s_1-\cdots-s_{\delta}} \ge {n-j-1\choose
          r-1}+1$, $\cdots$, or
          $d_{n-t_{\delta}-s_{\delta-\vartriangle}-\cdots-s_{\delta}} \ge {n-j-1\choose r-1}+\delta-\vartriangle$
or
          $d_{n-s_{\delta+1-\vartriangle}-\cdots-s_{\delta}} \ge {n-j-1\choose
          r-1}+\delta+1-\vartriangle $, $\cdots$, or $d_{n-s_{\delta}} \ge {n-j-1\choose
          r-1}+\delta$ or $d_n \ge {n-j-1\choose
          r-1}+\delta+1$  for  $j^\ast<j<{n\over 2}$,
          where $\vartriangle={n-j-1\choose  r-1}-{j-1\choose
          r-1}$.

        \item [$(3)$] $d_{n-t_1-\cdots-t_{\delta}-s_1-\cdots-s_{\delta}}\le
         {{n\over 2}-1\choose r-1}$, $d_{n-t_2-\cdots-t_{\delta}-s_2-\cdots-s_{\delta}}
         \le {{n\over 2}-1\choose r-1}+1$,$\cdots$, and $d_{n-t_{\delta}-s_{\delta}}
         \le {{n\over 2}-1\choose r-1}+\delta-1$ implies $d_n
         \ge {{n\over 2}-1\choose r-1}+\delta+1$ for even $n\ge 2r+2$,
where $(t_1+2t_2+\cdots+\delta t_{\delta})+(s_1+2s_2+\cdots+\delta
s_{\delta})=\delta r$,
        $\delta\le t_1+s_1\le \delta r$,
        $0\le t_2+s_2\le \lfloor{\delta\over 2}\rfloor(r-1)$, $\cdots$, $0\le t_{\delta}+s_{\delta}\le \lfloor{\delta\over {\delta}}\rfloor(r-1)$,
        $1\le t_1+t_2+\cdots+t_{\delta}\le \min\{\delta r-\delta, j\}$,     $1\le s_1+s_2+\cdots+s_{\delta}\le \min\{\delta r-\delta, n-j\}$,   $0\le t_2+\cdots+t_{\delta}\le \lfloor{\delta\over 2}\rfloor(r-1)$ if $\delta\ge 2$,
        $0\le t_3+\cdots+t_{\delta}\le \lfloor{\delta\over 3}\rfloor(r-1),\cdots$,
        $0\le t_{\delta-1}+t_{\delta}\le \lfloor{\delta\over \delta-1}\rfloor(r-1)$,
        $0\le s_2+\cdots+s_{\delta}\le \lfloor{\delta\over 2}\rfloor(r-1)$ if $\delta\ge 2$,
        $0\le s_3+\cdots+s_{\delta}\le \lfloor{\delta\over 3}\rfloor(r-1)$ if $\delta\ge 3$,$\cdots$,
        $0\le s_{\delta-1}+s_{\delta}\le \lfloor{\delta\over \delta-1}\rfloor(r-1)$ and
        $t_i( s_i)\ge 0$ for $1\le i\le \delta$.

        \end{itemize}
\end{theorem}

\pf If $d$ satisfies $(1)$ and $(3)$, then $d$ is forcibly maximally
edge-connected
 by Corollary \ref{28c}. Suppose that $d$ is not forcibly super edge-connected.
 Then there exists a non-super edge-connected $r$-uniform realization $H$,
 with $\lambda(H)=\delta$, consisting of two subhypergraphs $C_1$ and $C_2$
 joined by $\delta$ edges $e_1, \cdots, e_{\delta}$,
 where $V(C_2)=V(H)\setminus V(C_1)$ and $|V(C_1)|\le |V(C_2)|$.
 Let $|V(C_1)|=j$. Note that $r\le j \le \lfloor{n\over 2} \rfloor$.
Let $U_i(V_i)$ be the set of vertices in $V(H_1)(V(H_2))$, any one
of which belong to exactly $t_i(s_i)$ hyperedges of $e_1, e_2,
\cdots, e_{\delta}$, where $1\le i\le \delta$, then we know that
$(t_1+2t_2+\cdots+\delta t_{\delta})+(s_1+2s_2+\cdots+\delta
s_{\delta})=\delta r$, $\delta\le t_1+s_1\le \delta r$, $0\le
t_2+s_2\le \lfloor{\delta\over 2}\rfloor(r-1)$, $\cdots$, $0\le
t_{\delta}+s_{\delta}\le \lfloor{\delta\over {\delta}}\rfloor(r-1)$,
$1\le t_1+t_2+\cdots+t_{\delta}\le \min\{\delta r-\delta, j\}$,
$1\le s_1+s_2+\cdots+s_{\delta}\le \min\{\delta r-\delta, n-j\}$,
$0\le t_2+\cdots+t_{\delta}\le \lfloor{\delta\over 2}\rfloor(r-1),
\cdots$, $0\le t_{\delta-1}+t_{\delta}\le \lfloor{\delta\over
\delta-1}\rfloor(r-1)$, $0\le s_2+\cdots+s_{\delta}\le
\lfloor{\delta\over 2}\rfloor(r-1), \cdots$, $0\le
s_{\delta-1}+s_{\delta}\le \lfloor{\delta\over
\delta-1}\rfloor(r-1)$ and $t_i( s_i)\ge 0$ for $1\le i\le \delta$.

We can see that there are at most  $t_1+t_2+\cdots+t_{\delta}$
vertices in $C_1$ of degree larger than $j-1 \choose r-1$, at most
$t_2+\cdots+t_{\delta}$ vertices in $C_1$ of degree larger than
${j-1 \choose r-1}+1$,$\cdots$, and at most  $t_{\delta}$ vertices
in $C_1$ of degree larger than ${j-1 \choose r-1}+\delta-1$, while
the degree of each vertex in $C_1$ is at most ${j-1 \choose
r-1}+\delta$. Thus we have $d_{j}\le {j-1\choose r-1}+\delta$,
$d_{j-t_{\delta}}\le {j-1\choose r-1}+\delta-1$, $\cdots$,
$d_{j-t_2-\cdots-t_{\delta}}\le {j-1\choose r-1}+1$ and
$d_{j-t_1-t_2-\cdots-t_{\delta}}\le {j-1\choose r-1}$.

If $j<{n\over 2}$, by similar argument in the proof of Theorem
\ref{23t}, we can obtain some inequalities, which contradicts $(1)$
and $(2)$, respectively.

If $j={n\over 2}$ for even $n$, analogously,
 we have  $d_{n-t_1-t_2-\cdots-t_{\delta}-s_1-\cdots-s_{\delta}}\le {{n\over 2}-1\choose r-1}$,
  $d_{n-t_2-\cdots-t_{\delta}-s_2-\cdots-s_{\delta}}\le {{n\over 2}-1\choose r-1}+1,\cdots,$
   $d_{n-t_{\delta}-s_{\delta}}\le
   {{n\over 2}-1\choose r-1}+\delta-1$ and
   $d_n \le {{n\over 2}-1\choose r-1}+\delta$, contrary to $(3)$. \qed

When $t_1=(r-1)\delta, t_2=t_3=\cdots=t_{\delta}=0$, $s_1=
s_2=\cdots=s_{\delta-1}=0$ and $s_{\delta}=1$ in Theorem \ref{32t},
we obtain the following immediate consequence for $r$-uniform
hypergraphic sequence with minimum $\delta$ to be super
edge-connected, which generalize the relevant result of Liu et al.
[11] since $g(\delta, r)\ge (r-1)\delta$.

\begin{corollary}\label{33c}
    For integers $r\ge2$ and $\delta \ge 1$. Let $d=(d_1, d_2, \cdots, d_n)$ be an
    $r$-uniform hypergraphic sequence with $\delta=d_1\le d_2\le\cdots\le d_n$.
If $d_{j-(r-1)\delta}\le {j-1\choose r-1}$ and $d_j\le {j-1\choose
r-1}+\delta$ implies $d_{n-1}\ge{n-j-1\choose r-1}+\delta$ or
$d_{n}\ge{n-j-1\choose r-1}+\delta+1$ for $(r-1)\delta\le j \le
\lfloor{n\over 2}\rfloor$, then $d$ is forcibly super
edge-connected.
\end{corollary}

If $r=2$ in Corollary \ref{33c}, then we obtain the following result for a graphic sequence
with minimum element at least $1$ to be forcibly super edge-connected.

\begin{corollary} ([11])\label{34c}
    For integer $\delta \ge 1$. Let $d=(d_1, d_2, \cdots, d_n)$ be a graphic sequence with
     $\delta=d_1\le d_2\le\cdots\le d_n$.
If  $d_{j-\delta}\le {j-1}$ and $d_j\le j-1+\delta$ implies
$d_{n-1}\ge n-j-1+\delta$ or $d_{n}\ge n-j+\delta$ for $\delta\le
j\le \lfloor{n\over 2}\rfloor$, then $d$ is forcibly super
edge-connected.
\end{corollary}

\section{Another sufficient degree condition for forcibly k-edge-connected hypergraph}

In [14], Yin and Guo gave a sufficient degree condition for a
graphic sequence to be forcibly k-edge-connected. We shall extend
this result to hypergraphs. That is Theorem 4.2, which contains five
"Bondy-Chv$\acute{a}$tal type" conditions.

\begin{theorem}([14])\label{31t}
    Let $d=(d_1,d_2,\cdots,d_n)$ be a non-decreasing graphic
    sequence, and let $k\ge 2$ be an integer. If $d$ satisfies each
    of the following, then $d$ is forcibly $k$-edge-connected:

    \begin{itemize}
    \item [$(1)$] $d_1\ge k$,

    \item [$(2)$] for any integers $x,z,y$ and $j$ with $1\le x\le
    k-1$, $\lceil {k-1\over x}\rceil\le z\le k-x$,  $ z\le y\le k-1 $ and  $k+1\le j\le
    {n-z\over 2}$, we have that $d_{j-x}\le j-1$,
    $d_{j-x+1}\le \lfloor {k-1\over x}\rfloor+ j-1$ and
    $d_j\le z+j-1$ implies $d_{n-y}\ge n-j$ or $d_{n-y+1}\ge n-j+\lfloor {k-1\over
    y}\rfloor$ or $d_n\ge n-j+min\{x,k-y\}$,

    \item [$(3)$] for $n\ge 2k+2$ and $ {n-k+1\over 2}<j\le \lfloor{n \over
    2}\rfloor$, we have that $d_{j-1}\le j-1$ implies
$d_{n-k}\ge n-j$ or $d_{n-1}\ge n-j+1$ or $d_n\ge j-1+k$,

    \item [$(4)$] for $n\ge 2k+2$, $1\le y\le
    k-1$ and $n$ even, we have that $d_{{n\over 2}-k+1}\le {n\over 2}-1$
implies $d_{n-(k-1+y)}\ge {n\over 2}$ or
    $d_{n-y}\ge {n\over 2}+1$ or $d_{n-y+1}\ge max\{{n\over 2}+1,{n\over 2}+\lfloor {k-1\over
    y}\rfloor\}$ or $d_n\ge {n\over 2}+k-y$,

    \item [$(5)$] for $n\ge 2k+2$ and any integers $x,z,y$ and $j$
    with $2\le x\le k-2$, $\lceil {k-1\over x}\rceil\le z\le k-x$, $z\le y\le
    k-1$ and $ {n-z\over 2}<j\le \lfloor{n \over
    2}\rfloor$,

    \item [$(5.1)$] if $j<{n-z+min\{x,k-y\}\over 2}$, then we have that
    $d_{j-x}\le j-1$ implies
    $d_{n-(x+y)}\ge n-j$ or $d_{n-(x+y)+1}\ge max\{n-j,j+\lfloor{k-1\over x}\rfloor\}$
or $d_{n-y+1}\ge max\{z+j,n-j+\lfloor{k-1\over y}\rfloor\}$ or
$d_n\ge n-j+min\{x,k-y\}$,

    \item [$(5.2)$] if $j>{n-z+min\{x,k-y\}\over 2}$,
    then we have that $d_{j-x}\le j-1$ implies
    $d_{n-(x+y)}\ge n-j$ or $d_{n-(x+y)+1}\ge max\{n-j,j+\lfloor{k-1\over x}\rfloor\}$
or $d_{n-x+1}\ge max\{n-j+min\{x,k-y\},j+\lfloor{k-1\over
x}\rfloor\}$ or $d_n\ge z+j$,

    \item [$(5.3)$] if $j={n-z+min\{x,k-y\}\over 2}$, then we have that $d_{j-x}\le j-1$ implies
    $d_{n-(x+y)}\ge n-j$ or $d_{n-(x+y)+1}\ge max\{n-j,j+\lfloor{k-1\over x}\rfloor\}$
or $d_{n-(x+y)+2}\ge max\{n-j+\lfloor{k-1\over
y}\rfloor,j+\lfloor{k-1\over x}\rfloor\}$ or $d_n\ge z+j$.

    \end{itemize}
\end{theorem}

For integers $z\ge 1, r\ge 2$ and $r+1\le j\le \lfloor {n\over
2}\rfloor$, we define $j_{0}=j(z,r)(\ge r+1)$ to be the biggest
integer such that $z+{j_0-1\choose r-1}\le {n-j_0-1\choose r-1}$.
Clearly, $r+1\le j_{0}< {n \over 2}$. Otherwise, when $j_{0}={n
\over 2}$, the contradiction $z\le 0$ is derived.

\begin{theorem}\label{31t}
    Given integers $r\ge 2$ and $k\ge 2$. Let $d=(d_1, d_2, \cdots, d_n)$
    be an $r$-uniform hypergraphic sequence with $d_1\le d_2\le\cdots\le d_n$.
    If $d$ satisfies each of the following, then $d$ is forcibly k-edge-connected:

    \begin{itemize}
    \item [$(1)$] $d_1\ge k$,

    \item [$(2)$] for any integers $x,z,y$ and $j$ with $1\le x\le
    (k-1)(r-1)$, $\lceil {k-1\over q}\rceil\le z\le k-q$, $x\equiv
    q(mod(k-1))$, $1\le q\le k-1, z\le y\le (k-1)(r-1)$ and $r+1\le j\le
    j_{0}$, we have that $d_{j-x}\le {j-1\choose r-1}$,
    $d_{j-x+1}\le \lfloor {k-1\over q}\rfloor+ {j-1\choose r-1}$ and
    $d_j\le z+{j-1\choose r-1}$ implies $d_{n-y}\ge {n-j-1\choose
    r-1}+1$ or $d_{n-y+1}\ge {n-j-1\choose r-1}+\lfloor {k-1\over
    R}\rfloor+1$ or $d_n\ge {n-j-1\choose r-1}+min\{q,k-R\}+1$,
    where $y\equiv R(mod(k-1))$ and $1\le R\le k-1$.

    \item [$(3)$] for $n\ge 2r+2$ and any integers $x$ and $j$ with
$r-1\le x\le (k-1)(r-2)+1$ and $ j_{0}<j\le \lfloor{n \over
    2}\rfloor$,
we have that $d_{j-x}\le {j-1\choose r-1}$ implies
    $d_{n-(x+k-1)}\ge {n-j-1\choose r-1}+1$ or $d_{n-x}\ge {n-j-1\choose
    r-1}+2$ or $d_n\ge {j-1\choose r-1}+k$.

\item [$(4)$] for $n\ge 2r+2$ and any integers $x$, $y$ and $j$ with
$1\le x\le (k-1)(r-2)$, $k-1\le y\le (k-1)(r-1)$ and $ j_{0}<j\le
\lfloor{n \over
    2}\rfloor$,
we have that $d_{j-x}\le {j-1\choose r-1}$ implies
    $d_{n-(x+y)}\ge {n-j-1\choose r-1}+1$ or
    $d_{n-x}\ge {n-j-1\choose
    r-1}+2$ or $d_n\ge {j-1\choose r-1}+k$.

    \item [$(5)$] for $n\ge 2r+2$, $k-1\le x\le (k-1)(r-1)$, $k-1\le y\le
    (k-1)(r-1)$ and $n$ even, we have that $d_{{n\over 2}-x}\le {{n\over 2}-1\choose
    r-1}$ implies $d_{n-(x+y)}\ge {{n\over 2}-1\choose r-1}+1$ or
    $d_{n-y}\ge {{n\over 2}-1\choose
    r-1}+2$ or $d_{n-y+1}\ge max\{{{n\over 2}-1\choose r-1}+2,{{n\over 2}-1\choose r-1}+\lfloor {k-1\over
    R}\rfloor+1\}$ or $d_n\ge {{n\over 2}-1\choose r-1}+k-R+1$.

    \item [$(6)$] for $n\ge 2r+2$ and any integers $x,z,y$ and $j$
    with $2\le z\le k-2$, $1\le x\le (k-1)(r-1)$, $z\le y\le
    (k-1)(r-1)$ and $ j_{0}<j\le \lfloor{n \over
    2}\rfloor$,

    \item [$(6.1)$] if $z+{j-1\choose r-1}<{n-j-1\choose
    r-1}+min\{q,k-R\}$, then we have that
    $d_{j-x}\le {j-1\choose r-1}$ implies
    $d_{n-(x+y)}\ge {n-j-1\choose r-1}+1$ or $d_{n-(x+y)+1}\ge max\{{n-j-1\choose
    r-1},{j-1\choose r-1}+\lfloor{k-1\over q}\rfloor\}+1$ or $d_{n-y+1}\ge max\{{j-1\choose
    r-1}+z,{n-j-1\choose r-1}+\lfloor{k-1\over R}\rfloor\}+1$ or
     $d_n\ge {n-j-1\choose r-1}+min\{q,k-R\}+1$.

    \item [$(6.2)$] if $z+{j-1\choose r-1}>{n-j-1\choose
    r-1}+min\{q,k-R\}$, then we have that $d_{j-x}\le {j-1\choose r-1}$ implies
    $d_{n-(x+y)}\ge {n-j-1\choose r-1}+1$ or $d_{n-(x+y)+1}\ge max\{{n-j-1\choose
    r-1},{j-1\choose r-1}+\lfloor{k-1\over q}\rfloor\}+1$ or
$d_{n-x+1}\ge max\{{n-j-1\choose
    r-1}+min\{q,k-R\},{j-1\choose r-1}+\lfloor{k-1\over q}\rfloor\}+1$ or
$d_n\ge {j-1\choose r-1}+z+1$.

    \item [$(6.3)$] if $z+{j-1\choose r-1}={n-j-1\choose
    r-1}+min\{q,k-R\}$, then we have that $d_{j-x}\le {j-1\choose r-1}$ implies
    $d_{n-(x+y)}\ge {n-j-1\choose r-1}+1$
    or $d_{n-(x+y)+1}\ge max\{{n-j-1\choose r-1},
{j-1\choose r-1}+\lfloor{k-1\over q}\rfloor\}+1$ or
$d_{n-(x+y)+2}\ge max\{{j-1\choose r-1}+\lfloor{k-1\over
q}\rfloor,{n-j-1\choose r-1}+\lfloor{k-1\over R}\rfloor\}+1$ or
$d_n\ge {j-1\choose r-1}+z+1$.

    \end{itemize}
\end{theorem}

\pf Suppose $d$ satisfies (1)-(6) of Theorem 4.2, but is not
forcibly k-edge-connected. Then $d$ has a realization $H$ with an
edge-cut $E_0$ of $k-1$ hyperedges joining two subhypergraphs $H_1$
and $H_2$, with $|V(H_1)|=j$, $|V(H_2)|=n-j$ and $|V(H_1)|\le
|V(H_2)|$. Note that $r+1\le j\le \lfloor {n\over 2}\rfloor$ and
$n\ge 2r+2$, since $j\le r$ implies there exists a veretex of degree
less than 2, contradicting $d_1\ge k\ge 2$.

 By $E_H(V(H_1),V(H_2))$ we
denote the set of edges of $H$ which connect a vertex from $V(H_1)$
to one of $V(H_2)$. By $\partial _H V(H_1)$ we denote the set of
vertices of $H$ from $V(H_1)$ which have at least one neighbor
outside $V(H_1)$. Let $F=(X\cup Y,E_H(X,Y))$, $X=\partial _H V(H_1)$
and $Y=\partial _H V(H_2)$. Denote $x=|X|,y=|Y|$ and $z=max\{d_F
(v)|v\in X\}$. Then $1\le x\le
    (k-1)(r-1)$, $\lceil {k-1\over q}\rceil\le z\le k-q$
    and  $ z\le y\le (k-1)(r-1)$, where $x\equiv
    q(mod(k-1))$ and $1\le q\le k-1$. We consider two cases.

{\bf Case  1.} $ r+1\le j\le j_{0}$.

Let $z'=min\{d_F (v)|v\in X\}$, then $z'\le \lfloor {k-1\over
q}\rfloor$. At most $x$ vertices in $H_1$ can have degree larger
than ${j-1\choose r-1}$ and at most $x-1$ vertices in $H_1$ can have
degree larger than ${j-1\choose r-1}+z'$, and so $d_{j-x}\le
{j-1\choose r-1}$,
    $d_{j-x+1}\le \lfloor {k-1\over q}\rfloor+ {j-1\choose r-1}$.
    Also, no vertex in $H_1$ can have
degree larger than ${j-1\choose r-1}+z$, and so $d_j\le {j-1\choose
r-1}+z$. No vertex in $H_2 \backslash Y$ has degree at most
${n-j-1\choose r-1}$, and since ${j-1\choose r-1}+z\le {n-j-1\choose
r-1}$, we have $d_{n-y}\le {n-j-1\choose r-1}$.

Let $w=max\{d_F (v)|v\in Y\}$ and $w'=min\{d_F (v)|v\in Y\}$, then
$w\le min\{q,k-R\}$ and $w'\le \lfloor {k-1\over R}\rfloor$, where
$y\equiv R(mod(k-1))$ and $1\le R\le k-1$. Thus, at most $y-1$
vertices in $H_2$ can have degree larger than ${n-j-1\choose
r-1}+w'$ and each vertex in $H_2$ can has degree at most
${n-j-1\choose r-1}+w$. Thus $d_{n-y+1}\le {n-j-1\choose
r-1}+\lfloor {k-1\over R}\rfloor$ and $d_n\le {n-j-1\choose
r-1}+min\{q,k-R\}$. This contradicts (2).

{\bf Case  2.} $ j_{0}<j\le \lfloor{n \over 2}\rfloor$.

By the definition of $j_0$, we know that ${n-j-1\choose r-1}<
{j-1\choose r-1}+z $. Now three cases arise.

{\bf Subcase 2.1.}  $z=1$.

 Note that $k-1\le x\le (k-1)(r-1)$, $k-1\le y\le
    (k-1)(r-1)$. It follows clearly from ${j-1\choose r-1}\le
{n-j-1\choose r-1} $ and ${n-j-1\choose r-1}< {j-1\choose r-1}+1 $
that $n$ is even and $j={n\over 2}$.

At most $x$ vertices in $H_1$ can have degree larger than ${{n\over
2}-1\choose r-1}$ and at most $y$ vertices in $H_2$ can have degree
larger than ${{n\over 2}-1\choose r-1}$, so $d_{{n\over
2}-x}\le{{n\over 2}-1\choose r-1}$ and $d_{n-(x+y)}\le{{n\over
2}-1\choose r-1}$. No vertex in $H_1$ can have degree larger than
${{n\over 2}-1\choose r-1}+1$, and at most $y$ vertices in $H_2$ can
have degree larger than ${{n\over 2}-1\choose r-1}+1$ and at most
$y-1$ vertices in $H_2$ can have degree larger than ${{n\over
2}-1\choose r-1}+\lfloor {k-1\over R}\rfloor$, and so
$d_{n-y}\le{{n\over 2}-1\choose r-1}+1$ and $d_{n-y+1}\le
max\{{{n\over 2}-1\choose r-1}+1,{{n\over 2}-1\choose r-1}+\lfloor
{k-1\over R}\rfloor\}$. No vertex in $H_2$ can have degree larger
than ${{n\over 2}-1\choose r-1}+k-R$, and so $d_n\le {{n\over
2}-1\choose r-1}+k-R$. Thus (5) fails, a contradiction.

{\bf Subcase 2.2.} $z=k-1$.

{\bf Subsubcase 2.2.1.} $y=z$.

Note that $r-1\le x\le (r-2)(k-1)+1$. There are at most $x$ vertices
in $H_1$ can have degree larger than ${j-1\choose r-1}(\le
{n-j-1\choose r-1})$, at most $k-1$ vertices in $H_2$ can have
degree larger than ${n-j-1\choose r-1}$, and so $d_{j-x}\le
{j-1\choose r-1}$ and $d_{n-k+1-x}\le {n-j-1\choose r-1}$. No vertex
in $H_1$ can have degree larger than ${j-1\choose r-1}+k-1$, and no
vertex in $H_2$ can have degree larger than ${n-j-1\choose r-1}+1$.
Since $z+{j-1\choose r-1}>{n-j-1\choose r-1}$, $k-1+{j-1\choose
r-1}\ge {n-j-1\choose r-1}+1$, and so $d_{n-x}\le {n-j-1\choose
r-1}+1$ and $d_n\le {j-1\choose r-1}+k-1$. Thus (3) fails, a
contradiction.

{\bf Subsubcase 2.2.2.} $z<y\le (r-1)(k-1)$.

Note that $1\le x\le (r-2)(k-1)$. At most $x$ vertices in $H_1$ can
have degree larger than ${j-1\choose r-1}(\le {n-j-1\choose r-1})$,
at most $y(>k-1)$ vertices in $H_2$ can have degree larger than
${n-j-1\choose r-1}$, and so $d_{j-x}\le {j-1\choose r-1}$,
$d_{n-x-y}\le {n-j-1\choose r-1}$, $d_{n-x-y+1}\le {n-j-1\choose
r-1}+1$ and $d_{n-x}\le {n-j-1\choose r-1}+1$. No vertex in $H_1$
can have degree larger than ${j-1\choose r-1}+k-1$, and no vertex in
$H_2$ can have degree larger than ${n-j-1\choose r-1}+1(\le
{n-j-1\choose r-1}+k-1)$. So $d_n\le {j-1\choose r-1}+k-1$. Thus (4)
fails, a contradiction.

{\bf Subcase 2.3.}    $2\le z\le k-2.$

{\bf Subsubcase 2.3.1.}  $z+{j-1\choose r-1}<{n-j-1\choose
r-1}+min\{q,k-R\}$.

There are at most $x$ vertices in $H_1$ can have degree larger than
${j-1\choose r-1}(\le {n-j-1\choose r-1})$, and at most $y$ vertices
in $H_2$ can have degree larger than $ {n-j-1\choose r-1}$, and so
$d_{n-(x+y)}\le {n-j-1\choose r-1}$. At most $x-1$ vertices in $H_1$
can have degree larger than $max\{{n-j-1\choose
    r-1},{j-1\choose r-1}+\lfloor{k-1\over q}\rfloor\}$ and at most $y$ vertices in
$H_2$ can have degree larger than $max\{{n-j-1\choose
    r-1},{j-1\choose r-1}+\lfloor{k-1\over q}\rfloor\}$, and so
$d_{n-(x+y)+1}\le max\{{n-j-1\choose
    r-1},{j-1\choose r-1}+\lfloor{k-1\over q}\rfloor\}$.
    No vertex in $H_1$ can have
degree larger than $max\{{j-1\choose
    r-1}+z,{n-j-1\choose r-1}+\lfloor{k-1\over R}\rfloor\}$, and at most $y-1$ vertices in
$H_2$ can have degree larger than $max\{{j-1\choose
    r-1}+z,{n-j-1\choose r-1}+\lfloor{k-1\over R}\rfloor\}$, and so $d_{n-y+1} \le max\{{j-1\choose
    r-1}+z,{n-j-1\choose r-1}+\lfloor{k-1\over R}\rfloor\}$.
No vertex in $H_1$ can have degree larger than ${j-1\choose r-1}+z$
and no vertex in $H_2$ can have degree larger than ${n-j-1\choose
r-1}+min\{q,k-R\}$, and so $d_n\le {n-j-1\choose r-1}+min\{q,k-R\}$,
contrary to (6.1).

{\bf Subsubcase 2.3.2.} $z+{j-1\choose r-1}>{n-j-1\choose
r-1}+min\{q,k-R\}$.

At most $x$ vertices in $H_1$ can have degree larger than
${j-1\choose r-1}(\le {n-j-1\choose r-1})$, and at most $y$ vertices
in $H_2$ can have degree larger than $ {n-j-1\choose r-1}$, and so
$d_{n-(x+y)}\le {n-j-1\choose r-1}$. At most $x-1$ vertices in $H_1$
can have degree larger than $max\{{n-j-1\choose
    r-1},{j-1\choose r-1}+\lfloor{k-1\over q}\rfloor\}$ and at most $y$ vertices in
$H_2$ can have degree larger than $max\{{n-j-1\choose
    r-1},{j-1\choose r-1}+\lfloor{k-1\over q}\rfloor\}$, and so
$d_{n-(x+y)+1}\le max\{{n-j-1\choose
    r-1},{j-1\choose r-1}+\lfloor{k-1\over q}\rfloor\}$.
 No vertex in $H_2$ can have
degree larger than $max\{{n-j-1\choose
    r-1}+min\{q,k-R\},{j-1\choose r-1}+\lfloor{k-1\over q}\rfloor\}$, and at most $x-1$ vertices in
$H_1$ can have degree larger than $max\{{n-j-1\choose
    r-1}+min\{q,k-R\},{j-1\choose r-1}+\lfloor{k-1\over q}\rfloor\}$,  and so $d_{n-x+1}\le max\{{n-j-1\choose
    r-1}+min\{q,k-R\},{j-1\choose r-1}+\lfloor{k-1\over
    q}\rfloor\}$.
No vertex in $H_1$ can have degree larger than ${j-1\choose r-1}+z$
and no vertex in $H_2$ can have degree larger than ${n-j-1\choose
r-1}+min\{q,k-R\}$, and so $d_n\le {j-1\choose r-1}+z$, contrary to
(6.2).

{\bf Subsubcase 2.3.3.} $z+{j-1\choose r-1}={n-j-1\choose
r-1}+min\{q,k-R\}$.

There are at most $x$ vertices in $H_1$ can have degree larger than
${j-1\choose r-1}(\le {n-j-1\choose r-1})$, and at most $y$ vertices
in $H_2$ can have degree larger than $ {n-j-1\choose r-1}$, and so
$d_{n-(x+y)}\le {n-j-1\choose r-1}$. At most $x-1$ vertices in $H_1$
can have degree larger than $max\{{n-j-1\choose
    r-1},{j-1\choose r-1}+\lfloor{k-1\over q}\rfloor\}$ and at most $y$ vertices in
$H_2$ can have degree larger than $max\{{n-j-1\choose
    r-1},{j-1\choose r-1}+\lfloor{k-1\over q}\rfloor\}$, and so
$d_{n-(x+y)+1}\le max\{{n-j-1\choose
    r-1},{j-1\choose r-1}+\lfloor{k-1\over q}\rfloor\}$.
At most $x-1$ vertices in $H_1$ can have degree larger than
$max\{{n-j-1\choose
    r-1}+\lfloor{k-1\over R}\rfloor, {j-1\choose r-1}+\lfloor{k-1\over q}\rfloor\}$ and at most $y-1$ vertices in
$H_2$ can have degree larger than $max\{{n-j-1\choose
    r-1}+\lfloor{k-1\over R}\rfloor, {j-1\choose r-1}+\lfloor{k-1\over q}\rfloor\}$, and so
$d_{n-(x+y)+2}\le max\{{n-j-1\choose
    r-1}+\lfloor{k-1\over R}\rfloor, {j-1\choose r-1}+\lfloor{k-1\over q}\rfloor\}$.
No vertex in $H_1$ can have degree larger than ${j-1\choose
    r-1}+z$, and no vertex in
$H_2$ can have degree larger than ${n-j-1\choose r-1}+min\{q,k-R\}$,
     and so $d_n\le {j-1\choose r-1}+z$,
contrary to (6.3). \qed

When $k=2$ in Theorem 4.2, we obtain $(z,q,R)=(1,1,1)$ and the
following conclusion.

\begin{corollary} \label{34c}
    Given integer $r \ge 2$. Let $d=(d_1, d_2, \cdots, d_n)$ be an
    r-uniform hypergraphic sequence with $d_1\le d_2\le\cdots\le d_n$.
If $d$ satisfies the following conditions, then $d$ is forcibly
2-edge-connected:

 \begin{itemize}
    \item [$(1)$] $d_1\ge 2$,

    \item [$(2)$] for any integers $x,y$ and $j$ with $1\le x\le r-1,1\le y\le
    r-1$ and $r+1\le j< {n\over 2}$, we have that $d_{j-x}\le {j-1\choose
    r-1}$ and $d_j\le {j-1\choose r-1}+1$ implies $d_{n-y}\ge {n-j-1\choose
    r-1}+1$ or $d_n\ge {n-j-1\choose r-1}+2$,

    \item [$(3)$]  for $n\ge 2r+2$ and $n$ even, we have that
     $d_{{n\over 2}}\le {{n\over 2}-1\choose r-1}$ implies
$d_{n-2}\ge {{n\over 2}-1\choose r-1}+1$ or $d_n\ge {{n\over
2}-1\choose r-1}+2$.

    \end{itemize}

\end{corollary}

If $k=3$, then from condition (2) of Theorem 4.2, we obtain
$(q,R,z)=(1,2,2),(2,2,1)$. Thus Theorem 4.2 reduces to Corollary
4.4.

\begin{corollary} \label{34c}
    Given integer $r \ge 2$. Let $d=(d_1, d_2, \cdots, d_n)$ be an
    r-uniform hypergraphic sequence with $d_1\le d_2\le\cdots\le d_n$.
If $d$ satisfies the following conditions, then $d$ is forcibly
3-edge-connected:

 \begin{itemize}
    \item [$(1)$] $d_1\ge 3$,

    \item [$(2)$] for any integers $x,y$ and $j$ with $1\le x\le 2(r-1),2\le y\le
    2(r-1)$ and $r+1\le j\le j_0$, we have that $d_{j-x}\le {j-1\choose
    r-1}$ and $d_j\le {j-1\choose r-1}+2$ implies $d_{n-y}\ge {n-j-1\choose
    r-1}+1$ or $d_n\ge {n-j-1\choose r-1}+2$,

\item [$(3)$] for any integers $x,y$ and $j$ with $1\le x\le 2(r-1),1\le y\le
    2(r-1)$ and $r+1\le j\le j_0$, we have that $d_{j-x}\le {j-1\choose
    r-1}$ and $d_j\le {j-1\choose r-1}+1$ implies $d_{n-y}\ge {n-j-1\choose
    r-1}+1$ or $d_n\ge {n-j-1\choose r-1}+2$,

\item [$(4)$] for  $n\ge 2r+2$ and   $j_0< j\le \lfloor {n\over 2}\rfloor$,
we have that $d_j\le {j-1\choose r-1}$ implies $d_{n-3}\ge
{n-j-1\choose
    r-1}+1$ or $d_{n-1}\ge
{n-j-1\choose
    r-1}+2$ or $d_n\ge {n-j-1\choose r-1}+3$,

    \item [$(5)$]  for $n\ge 2r+2$ and $n$ even,
$2\le x\le 2(r-1),2\le y\le
    2(r-1)$,
we have that
     $d_{{n\over 2}-x}\le {{n\over 2}-1\choose r-1}$ implies
$d_{n-(x+y)}\ge {{n\over 2}-1\choose r-1}+1$ or $d_n\ge {{n\over
2}-1\choose r-1}+2$.

    \end{itemize}

\end{corollary}

\end{document}